\documentclass{article}
\usepackage{amsmath,amssymb,amsthm}

\newtheorem{model}{Model}
\newtheorem*{cuh}{Computable Universe Hypothesis}

\begin{document}

\title{Is Turing's Thesis the Consequence of a More General Physical Principle?}

\author{Matthew~P. Szudzik}

\date{2 April 2012}

\maketitle

\begin{abstract}
We discuss historical attempts to formulate a physical hypothesis from which Turing's thesis may be derived, and also discuss some related attempts to establish the computability of mathematical models in physics.  We show that these attempts are all related to a single, unified hypothesis.
\end{abstract}

\section{Introduction} \label{s:intro}

Alan Turing~\cite{aT36} proposed the concept of a computer---that is, the concept of a mechanical device that can be programmed to perform any conceivable calculation---after studying the processes that humans use to perform calculations.  In particular, he claimed that any function of non-negative integers which can be effectively calculated by humans is a function that can be calculated by a Turing machine.  This claim, known as \emph{Turing's thesis}, is an empirical principle that has withstood many tests to its validity.\footnote{Indeed, every time a software developer attempts to write a computer program to implement an explicitly-described procedure, Turing's thesis is tested.}  But why does Turing's thesis seem to be true?  Can it be derived, for example, from a principle of contemporary sociology, from a principle of human biology, or perhaps from a principle of fundamental physics?

We will be concerned with the last of these questions, namely the question of whether Turing's thesis is the consequence of a physical principle.  But before discussing some of the historically important attempts to answer this question, let us introduce the following formalism.  Given a physical system, define a \emph{deterministic physical model} for the system to be
\begin{enumerate}
\item a set $S$ of \emph{states}
\end{enumerate}
together with
\begin{enumerate}
\setcounter{enumi}{1}
\item a set $A$ of functions from $S$ to the real numbers, and
\item for each non-negative real number $t$, a function $\varPhi_t$ from $S$ to $S$.
\end{enumerate}
If the system begins in state $s$, then we understand $\varPhi_t\left(s\right)$ to be the state of the system after $t$ units of time.  Furthermore, we identify each member $\alpha$ of $A$ with an \emph{observable quantity} of the system,\footnote{To ensure that the predictions of the model are unambiguous, it is common practice in physics to define each observable quantity operationally~\cite{pB27}.} and we consider $\alpha\left(s\right)$ to be the \emph{value predicted} for the observable quantity when the system is in state $s$.  For example, the following is a deterministic physical model for a particle that is moving uniformly along a straight line with a velocity of $3$ meters per second.
\begin{model} \label{m:detline}
Let $S$ be the set of all real numbers and define $\varPhi_t\left(x\right)=x+3t$ for all states $x$ in $S$, where $t$ is the time measured in seconds.  The position of the particle on the line, measured in meters, is given by the function $\alpha\left(x\right)=x$.
\end{model}
\noindent
We say that a deterministic physical model is \emph{faithful} if and only if there is a state $s_0$ such that, for each time $t$, the actual values of the observable quantities at that time agree with the values that are predicted when the system is in state $\varPhi_t\left(s_0\right)$.  Deterministic physical models are commonly studied in the theory of dynamical systems.

Now, the first notable attempt to derive Turing's thesis from a principle of physics appears to have been made by Robert Rosen~\cite{rR62}.  Rosen hypothesized\footnote{Rosen did not fully formalize his hypotheses.  The hypotheses given here are reasonable formal interpretations of Rosen's ``Hypothesis I''~\cite{rR62}.} that every physical system has a faithful deterministic physical model where
\begin{enumerate}
\item the set $S$ is the set of all non-negative integers,
\item \label{i:observable} each $\alpha$ in $A$ is a total recursive function,\footnote{We will use terminology from the theory of recursive functions~\cite{hR67} throughout this paper.  It follows from Turing's work~\cite{aT37} that the set of all functions from non-negative integers to non-negative integers that are computable by Turing machines is identical to the set of recursive functions.} and
\item when restricted to non-negative integer times $t$, $\varPhi_t\left(s\right)$ is a total recursive function of $s$ and $t$.
\end{enumerate}
Note that since recursive functions are functions from non-negative integers to non-negative integers, Condition~\ref{i:observable} implies that the values of all observable quantities are non-negative integers.  This hypothesis is justified by the fact that actual physical measurements have only finitely many digits of precision.  For example, a distance measured with a meterstick is a non-negative integer multiple $m$ of the length of the smallest division on the meterstick.  Therefore, without loss of generality, we can consider $m$ to be the value of the measurement.  Nevertheless, Condition~\ref{i:observable} does not forbid time from being measured with a non-negative real number, since time is not regarded as an observable quantity in deterministic physical models.  For example, it is consistent with Rosen's hypotheses to specify that
\begin{equation}
\varPhi_t\left(s\right)=\varPhi_{\left\lfloor t\right\rfloor}\left(s\right)
\end{equation}
for all non-negative real numbers $t$ and for all states $s$, where $\left\lfloor t\right\rfloor$ denotes the largest integer less than or equal to $t$.

Now, Turing's thesis can be derived from Rosen's hypotheses as follows.  First, note that by definition, if a function $\psi$ is \emph{effectively calculable} then there is a physical system that can reliably be used to calculate $\psi\left(x\right)$ for every non-negative integer $x$.  This means that there must be a system where the input can be observed, where the output can be observed, and where it is possible to observe that the system is finished with the calculation.  In the context of deterministic physical models this means that the system has observable quantities $\alpha$, $\beta$, and $\gamma$, respectively, such that for each non-negative integer $x$
\begin{enumerate}
\item there exists a state $s$ and a time $u$ such that $\alpha\left(s\right)=x$ and $\gamma\left(\varPhi_t\left(s\right)\right)=1$ for all $t\geq u$, and
\item if $\alpha\left(s\right)=x$ for any state $s$ then $\beta\left(\varPhi_t\left(s\right)\right)=\psi\left(x\right)$ where $t$ is any time such that $\gamma\left(\varPhi_t\left(s\right)\right)=1$.
\end{enumerate}
It then immediately follows from Rosen's hypotheses that every effectively calculable function is recursive, whether it is calculated by a human being or by any other physical system.  That is, Turing's thesis is a consequence of Rosen's hypotheses.

Informally, Rosen's hypotheses can be understood as stating that the universe is discrete, deterministic, and computable.  It is important then to ask whether we are living in a discrete, deterministic, computable universe.  Unfortunately, according to the currently-understood laws of physics, the answer appears to be ``No.''  The universe, as we currently understand it, does not seem to satisfy Rosen's hypotheses.  Nevertheless, it may be possible to reformulate the currently-understood physical laws so that Rosen's hypotheses are satisfied.  This sort of reformulation was first attempted by Konrad Zuse~\cite{kZ67,kZ69} in the 1960's.  Since then, increasingly sophisticated attempts have been made by Edward Fredkin~\cite{eF90} and Stephen Wolfram~\cite{sW02}.  In contrast, Roger Penrose~\cite{rP89,rP94} has speculated that the universe might not be computable, but efforts to find experimental evidence for this assertion have not succeeded.  There is, in fact, very little that can currently be said.  The question of whether Turing's thesis is the consequence of a valid physical principle is too difficult to be answered conclusively at this time.

\section{Physical Models}

Although deterministic physical models are widely used in the study of dynamical systems and have important real-world applications, they are not the only sort of model that is useful in physics.  The sorts of models used in quantum electrodynamics~\cite{rF85}, for example, are necessarily non-deterministic.  And it is difficult to reconcile the way that time is modeled in general relativity~\cite{kG49} with the way that it is treated as a linear quantity external to the states in a deterministic physical model.  Anyone who has attempted to reformulate the laws of physics so that Rosen's hypotheses are satisfied has had to grapple with these obstacles, and no one has had complete success.

For this reason, we have proposed~\cite{mS12} a more general sort of model which we simply call a \emph{physical model}.  A physical model for a system is a set $S$ of states together with a set $A$ of functions from $S$ to the real numbers.  Each member $\alpha$ of $A$ is identified with an observable quantity of the system,\footnote{As before, we require that each observable quantity be defined operationally.} and $\alpha\left(s\right)$ is the value predicted for that observable quantity when the system is in state $s$.  Deterministic physical models are special sorts of physical models.  For example, the deterministic physical model that was described in the introduction (Model~\ref{m:detline}) can be expressed as the following physical model.
\begin{model} \label{m:line}
Let $S$ be the set of all triples $\left(x,x_0,t\right)$ of real numbers such that $t>0$ and $x=x_0+3t$.  The position of the particle on the line, measured in meters, is given by the function $\alpha\left(x,x_0,t\right)=x$.  The initial position of the particle (for example, recorded in the observer's notebook) is given by the function $\beta\left(x,x_0,t\right)=x_0$.  The time, measured in seconds, is given by the function $\gamma\left(x,x_0,t\right)=t$.
\end{model}
\noindent
In contrast to Model~\ref{m:detline}, time is treated as an observable quantity in this model, as is the initial position.  The inclusion of these observable quantities in the model can be justified physically by noting that if a researcher were to test Model~\ref{m:detline} in a laboratory experiment, he would be required to measure the initial position $x_0$ of the particle and the position $x$ of the particle at some later time $t$.  That is, besides the position $x$, observations of both the time and the initial position are fundamental to the system.

Physical models can also be used to describe non-deterministic systems.  For example, suppose that one atom of the radioactive isotope nitrogen-$13$ is placed inside a detector at time $t=0$.  We say that the detector has status $1$ at time $t$ if the decay of the isotope was detected at any earlier time, and we say that the detector has status $0$ otherwise.  We use a non-negative integer to represent the history of the detector.  In particular, if $b_i$ is the status of the detector at time $t=i$, then the history of the detector at time $t=n$ is $h=\sum_{i=1}^nb_i2^{n-i}$.  Note that when $h$ is written as an $n$-bit binary number, the $i$th bit (counting from left to right) is $b_i$.  For example, if the isotope decays sometime between $t=2$ and $t=3$, then the history of the detector at time $t=5$ is $7$ because the binary representation of $7$ is $\left(00111\right)_2$.  Now, if the history of the detector is displayed on a computer screen, then the following is a physical model for the history of the detector as observed by a researcher looking at the screen.
\begin{model} \label{m:radio}
Let $S$ be the set of all triples $\left(t,2^d-1,k\right)$ where $d$, $t$, and $k$ are non-negative integers such that $d\leq t$, $t\ne0$, and $2k\leq2^d-1$.  The history of the detector is given by the function $\alpha\left(t,h,k\right)=h$, and the time, measured in units of the half-life of the isotope, is given by the function $\beta\left(t,h,k\right)=t$.
\end{model}

In a deterministic model such as Model~\ref{m:line}, the state of the system at a particular time is uniquely determined from its initial state.  But this is not the case for Model~\ref{m:radio}, which is a non-deterministic physical model.  For example, if we are given as an initial condition that the detector had status $0$ at time $t=1$, then there are four different states at time $t=3$ which are consistent with that initial condition.
\begin{center}
$\left(3,\left(000\right)_2,0\right)$\qquad$\left(3,\left(001\right)_2,0\right)$\qquad$\left(3,\left(011\right)_2,0\right)$\qquad$\left(3,\left(011\right)_2,1\right)$
\end{center}
In addition, Model~\ref{m:radio} has the property that if these states are considered to be equally likely, then the corresponding probability of observing a given history at time $t=3$ matches the probability predicted by the conventional theory of radioactive decay.  For example, since the history $\left(011\right)_2$ is observed in half of the four states, there is a $\frac{1}{2}$ probability that the detector's history will be $\left(011\right)_2$ at time $t=3$ if the detector had status $0$ at time $t=1$.  See~\cite{mS12} for a more complex example of a non-deterministic physical model that involves incompatible quantum measurements.

\section{Computable Physical Models}

Define a \emph{computable physical model} to be a physical model where $S$ is a recursive set of non-negative integers and where each $\alpha$ in $A$ is a total recursive function.\footnote{See the isomorphism theorems in~\cite{mS12} for other, equivalent definitions of a computable physical model.}  Of course, deterministic physical models that satisfy Rosen's hypotheses can all be expressed as computable physical models.  Non-deterministic models, such as the model for radioactive decay (Model~\ref{m:radio}), can also be expressed as computable physical models.  Now, we assert the following hypothesis.
\begin{cuh}
The laws of physics can be expressed as a computable physical model.
\end{cuh}
\noindent
To show that Turing's thesis is a consequence of the computable universe hypothesis, first recall that if a function $\psi$ is \emph{effectively calculable} then there is a physical system that can reliably be used to calculate $\psi\left(x\right)$ for every non-negative integer $x$.  This means that in some state of the universe it must be possible to observe a record of the the function's input, to observe the function's output, and to observe that the calculation has finished.  In the context of a physical model for the universe, this means that there are observable quantities $\alpha$, $\beta$, and $\gamma$, respectively, such that for each non-negative integer $x$
\begin{enumerate}
\item there exists a state $s$ such that $\alpha\left(s\right)=x$ and $\gamma\left(s\right)=1$, and
\item if $\alpha\left(s\right)=x$ and $\gamma\left(s\right)=1$ for any state $s$, then $\beta\left(s\right)=\psi\left(x\right)$.
\end{enumerate}
It then immediately follows from the computable universe hypothesis that every effectively calculable function is recursive, whether it is calculated by a human being or by any other system that is governed by physical law.  Hence, Turing's thesis is a consequence of the computable universe hypothesis.

Of course, as was the case with Rosen's hypotheses, it is not known whether the computable universe hypothesis is true.  But since the computable universe hypothesis allows for non-determinism and for more complex temporal relationships, it may be somewhat easier to reformulate the currently-understood physical laws so as to satisfy this less-restrictive hypothesis.

We conclude this section by considering the following example that helps to clarify certain features of our definition of effective calculability in the context of non-deterministic computable physical models.  First, let $\left\langle x,y\right\rangle$ denote a non-negative integer that encodes the pair $\left(x,y\right)$ of non-negative integers.  For example, using Cantor's pairing function, we could define
\begin{equation}
\left\langle x,y\right\rangle=\frac{1}{2}\left(x^2+2x y+y^2+3x+y\right)\enspace.
\end{equation}
Triples $\left(x,y,z\right)$ of non-negative integers can then be encoded as $\left\langle\left\langle x,y\right\rangle,z\right\rangle$.  We will use $\left\langle x,y,z\right\rangle$ as an abbreviation for $\left\langle\left\langle x,y\right\rangle,z\right\rangle$.  Also define $x\bmod2$ to be the rightmost bit in the binary representation of the non-negative integer $x$.  Now consider the following computable physical model for the history of a `noisy' detector.
\begin{model} \label{m:noise}
Let $S$ be the set of all pairs $\left\langle t,h\right\rangle$ where $t$ and $h$ are non-negative integers such that $t\ne0$ and $h\leq2^t-1$.  The time is given by the function $\alpha\left\langle t,h\right\rangle=t$, the current status of the detector is given by the function $\beta\left\langle t,h\right\rangle=h\bmod2$, a trivial observable quantity is given by the function $\gamma\left\langle t,h\right\rangle=1$, and the history of the detector is given by the function $\delta\left\langle t,h\right\rangle=h$.
\end{model}
\noindent
Although superficially similar to Model~\ref{m:radio}, all histories are possible in Model~\ref{m:noise}.  That is, there are no restrictions on the status of the detector.  At any point in time and regardless of the detector's past history, the status of the detector can be $1$ or $0$.  Now imagine forming a tree by taking each state $\left\langle t,h\right\rangle$ of Model~\ref{m:noise} as a node of the tree, and by taking the children of this node to be those states of the form $\left\langle t+1,h^\prime\right\rangle$ where the first $t$ bits of $h^\prime$ agree with the first $t$ bits of $h$.  We call this the \emph{tree of histories} for Model~\ref{m:noise}, and every branch on this tree corresponds to an alternate sequence of histories for the detector.  Moreover, we can think of the status of the detector as defining a function along each branch of the tree.  For each state $s$ on the branch, the input of the function is $\alpha\left(s\right)$, the corresponding output is $\beta\left(s\right)$, and the fact that the calculation has finished is indicated by $\gamma\left(s\right)$.  But since all histories are possible for this detector, every possible function $\psi$ from the positive integers to the set $\left\{0,1\right\}$ is calculated by the detector along some branch, including functions $\psi$ which are not recursive (that is, not computable by a Turing machine).  Therefore, Model~\ref{m:noise} is an example of a computable physical model where non-recursive functions may be calculated along certain branches of the tree of histories.  It is important to note, though, that according to our definition of effective calculability, these non-recursive functions are not effectively calculable, since for each state $s$ of Model~\ref{m:noise} there is another state $s^\prime$ such that $\alpha\left(s\right)=\alpha\left(s^\prime\right)$ and $\gamma\left(s\right)=\gamma\left(s^\prime\right)$, but $\beta\left(s\right)\ne\beta\left(s^\prime\right)$.  In other words, the detector cannot reliably be used to calculate $\psi$ because the detector is behaving non-deterministically.

\section{Continuous Models}

Although Turing's thesis holds in all discrete, deterministic, computable universes, Turing himself did not believe that the universe is discrete.  In particular, Turing~\cite{aT50} stated that
\begin{quote}
digital computers~\ldots~may be classified amongst the `discrete state machines'.  These are the machines which move by sudden jumps or clicks from one quite definite state to another.  These states are sufficiently different for the possibility of confusion between them to be ignored.  Strictly speaking there are no such machines.  Everything really moves continuously.
\end{quote}
But discreteness is not a prerequisite for computability.  In fact, Georg Kreisel~\cite{gK74} has hypothesized that the universe may be continuous and computable.

Before we turn to a discussion of Kreisel's hypothesis, define a \emph{functional physical model} for a system to be a set $D$ of finitely many real-valued functions, each of which takes $k$ real numbers as input, for some non-negative integer $k$.  We identify each member $\delta$ of $D$ with an observable quantity of the system, and also identify each of the $k$ inputs with an observable quantity.  The observable quantities that are identified with the inputs are said to be the \emph{given quantities} for the model,\footnote{To ensure that each function $\delta$ in $D$ is defined for all real number inputs, the scale of each given quantity should be adjusted so that it ranges over all real numbers.  For example, if a given quantity is a temperature and if $\delta$ is undefined for temperatures below $0$ degrees Kelvin, then a logarithmic temperature scale should be used instead of the Kelvin scale.} and the observable quantities that are identified with the members of $D$ are said to be the \emph{predicted quantities} of the model.  We say that a functional physical model is \emph{faithful} if and only if the values of the predicted quantities in the model match the values that are actually measured whenever the values of the given quantities in the model match the values that are actually measured.\footnote{In practice, of course, measurement errors will prevent us from knowing these values exactly.}

For example, the deterministic physical model that was described in the introduction (Model~\ref{m:detline}) can be expressed as the following functional physical model.
\begin{model} \label{m:func}
Let the position of the particle on the line, measured in meters, be given by the function $\delta\left(x_0,t\right)=x_0+3t$.  The initial position of the particle is $x_0$.  The time, measured in seconds, is $t$.
\end{model}
\noindent
In this case, the position of the particle is the only predicted quantity.  The particle's initial position and the time are the given quantities.

Now, note that every integer $i$ can be encoded as a non-negative integer $\zeta\left(i\right)$, where $\zeta\left(i\right)=2i$ if $i\geq0$ and where $\zeta\left(i\right)=-2i-1$ if $i<0$.  Each rational number $\frac{a}{b}$ in lowest-terms with $b>0$ can be encoded as a non-negative integer $\rho\left(\frac{a}{b}\right)$ where
\begin{equation}
\rho\left(\frac{a}{b}\right)=\zeta\left(\left(\mathrm{sgn}\:a\right)2^{\zeta\left(a_1-b_1\right)} 3^{\zeta\left(a_2-b_2\right)} 5^{\zeta\left(a_3-b_3\right)} 7^{\zeta\left(a_4-b_4\right)} 11^{\zeta\left(a_5-b_5\right)}\cdots\right)
\end{equation}
and where
\begin{equation}
a=\left(\mathrm{sgn}\:a\right)2^{a_1}3^{a_2}5^{a_3}7^{a_4}11^{a_5}\cdots
\end{equation}
and
\begin{equation}
b=2^{b_1}3^{b_2}5^{b_3}7^{b_4}11^{b_5}\cdots
\end{equation}
are the prime factorizations of $a$ and $b$, respectively.  For each pair of rational numbers $q$ and $r$, define $\left(q\,;r\right)$ to be the non-negative integer $\left\langle\rho\left(q\right),\rho\left(r\right)\right\rangle$.  We will use $\left(q\,;r\right)$ to represent the open interval with endpoints $q$ and $r$.

Next, let $\mathbb{N}$ denote the set of non-negative integers and let $\mathbb{R}$ denote the set of real numbers.  Note that every real number $x$ can be represented as a nested sequence of open intervals whose intersection is $x$.  We say that a function $\phi:\mathbb{N}\to\mathbb{N}$ is a \emph{nested oracle} for $x\in\mathbb{R}$ if and only if $\phi\left(0\right)=\left(a_0\,;b_0\right)$, $\phi\left(1\right)=\left(a_1\,;b_1\right)$, $\phi\left(2\right)=\left(a_2\,;b_2\right)$,~\ldots~is a sequence of nested intervals whose intersection is $x$.  A function $\delta:\mathbb{R}\to\mathbb{R}$ is said to be \emph{computable} (in the sense defined by Lacombe~\cite{dL55a,dL55b,dL55c}) if and only if there is a total recursive function $\xi$ such that if $\phi$ is a nested oracle for $x\in\mathbb{R}$, then $\lambda m\left[\,\xi\left(\phi\left(m\right)\right)\right]$ is a nested oracle for $\delta(x)$.  This definition can naturally be extended to functions from $\mathbb{R}^k$ to $\mathbb{R}$, for any $k\in\mathbb{N}$.  Finally, we say that a functional physical model is \emph{computable} if and only if every member of the set $D$ is computable in the sense that we have just described.

Let us now return to a discussion of Kreisel's hypothesis.  Kreisel hypothesized that every faithful functional physical model is computable.  For example, Model~\ref{m:func} is computable because in that model the function $\delta:\mathbb{R}^2\to\mathbb{R}$ is computed by the total recursive function $\xi$ that is defined so that
\begin{equation}
\xi\left(\left(a\,;b\right),\left(c\,;d\right)\right)=\left(a+3c\,;b+3d\right)
\end{equation}
for all rational numbers $a$, $b$, $c$, and $d$.  Moreover, every computable functional physical model can be expressed as a computable physical model.  For example, Model~\ref{m:func} can be expressed as the following computable physical model.
\begin{model} \label{m:contline}
Let $S$ be the set of all triples $\left\langle \left(a\,;b\right),\left(c\,;d\right),\left(a+3c\,;b+3d\right)\right\rangle$ where $a$, $b$, $c$, and $d$ are rational numbers such that $a<b$ and $c<d$.  The initial position of the particle on the line, represented as a range of positions measured in meters, is given by the function $\alpha\left\langle \left(a\,;b\right),\left(c\,;d\right),\left(e\,;f\right)\right\rangle=\left(a\,;b\right)$.  The time interval, measured in seconds, is given by the function $\beta\left\langle \left(a\,;b\right),\left(c\,;d\right),\left(e\,;f\right)\right\rangle=\left(c\,;d\right)$.  The position of the particle, represented as a range of positions measured in meters, is given by the function $\gamma\left\langle \left(a\,;b\right),\left(c\,;d\right),\left(e\,;f\right)\right\rangle=\left(e\,;f\right)$.
\end{model}
\noindent
Note that the states of Model~\ref{m:contline} are all of the form
\begin{equation}
\left\langle\left(a\,;b\right),\left(c\,;d\right),\xi\left(\left(a\,;b\right),\left(c\,;d\right)\right)\right\rangle\enspace,
\end{equation}
and this guarantees that if $\phi$ is a nested oracle for some initial position $x_0$ and if $\psi$ is a nested oracle for some time $t$, then there is a unique sequence of states $s_0$, $s_1$, $s_2$,~\ldots~in $S$ such that $\alpha\left(s_m\right)=\phi\left(m\right)$ and $\beta\left(s_m\right)=\psi\left(m\right)$ for all $m\in\mathbb{N}$.  Therefore,
\begin{equation}
\gamma\left(s_m\right)=\xi\left(\phi\left(m\right),\psi\left(m\right)\right)
\end{equation}
for all $m\in\mathbb{N}$, and $\lambda m\left[\gamma\left(s_m\right)\right]$ is a nested oracle for the position $\delta\left(x_0,t\right)$ that is predicted by Model~\ref{m:func}.  Thus, there is a direct correspondence between the functional physical model (Model~\ref{m:func}) and the computable physical model (Model~\ref{m:contline}).  See~\cite{mS12} for more information regarding this correspondence.

In summary, the computable physical models comprise a very general class of models, capable of expressing discrete deterministic models such as those studied by Rosen, non-deterministic models such as the model for radioactive decay, and continuous models such as those studied by Kreisel.  Furthermore, Turing's thesis is a consequence of the hypothesis that the laws of physics can be expressed as a computable physical model.  It is very tempting, therefore, to wonder whether this hypothesis might be true.

\bibliographystyle{amsplain}
\bibliography{IsTuringsThesis}

\end{document}